\documentclass[11pt]{amsart}

\usepackage{amsmath} 
\usepackage{amssymb}

\newtheorem{theorem}{Theorem}
\newtheorem{proposition}[theorem]{Proposition}

\theoremstyle{definition}

\theoremstyle{remark}
\newtheorem{conrem}[theorem]{Concluding Remarks}
\newtheorem{notation}[theorem]{Notation}

\newcommand{\metricD}{{\mathfrak d}}
\newcommand{\length}{{\rm length}}
\newcommand{\lh}{{\rm lg}}
\newcommand{\Aut}{{\rm Aut}}
\newcommand{\conc}{{}^\frown\!}

\title[Countable structure does not have$\ldots$]{Countable structure does
not have a free uncountable automorphism group} 
\author{Saharon Shelah}
\address{Institute of Mathematics\\
 The Hebrew University of Jerusalem\\
 Jerusalem 91904, Israel\\
 and Institute Mittag-Leffler\\
 The Royal Swedish Academy of Sciences\\
 Auravagen 17, S182 62\\
 Jursholm, Sweden\\
 and  Department of Mathematics\\
 Rutgers University\\
 New Brunswick, NJ 08854, USA}
\email{shelah@math.huji.ac.il}
\urladdr{http://www.math.rutgers.edu/$\sim$shelah}
\thanks{This research was partially supported by the Israel Science
 Foundation founded by the Israel Academy of Sciences and
 Humanities. Publication 744}   
\subjclass{}
\keywords{} 

\begin{document}
\maketitle

This is a well known problem in group theory and we thank Simon Thomas for
telling us about it. Solecki \cite{So99} proved that the group of
automorphisms of a countable structure cannot be an uncountable free abelian
group.  See more in Just, Shelah and Thomas \cite{JShT:654} where, as a
byproduct, we can say something on uncountable structures. 
  
Here, we prove 
  
\begin{theorem}
\label{Thm1}
If ${\mathbb A}$ is a countable model, then $\Aut({\mathbb A})$ cannot be a
free uncountable group.
\end{theorem}

The proof follows from the following two claims, one establishing a property
of $G$ and the other proving that free groups does not have it.  
  
\begin{notation}
\begin{enumerate}
\item Let $\omega$ denote the set of natural numbers, and let $x<\omega$
mean ``$x$ is a natural number''.  
\item Let $a,b,c,d$ denote members of $G$ (the group).
\item Let ${\mathbf d}$ denote the $\omega$--sequence $\langle d_{n}:n<\omega
\rangle$, and similarly in other cases.  
\item Let $k,\ell,m,n,i,j,r,s,t$ denote natural numbers (and so also
elements of the structure ${\mathbb A}$, which we assume is the set of
natural number for notational simplicity).
\end{enumerate}
\end{notation}

\begin{proposition}
\label{prop3}
Assume ${\mathbb A}$ a countable structure with automorphism group $G$, and
for notational simplicity its set of elements is $\omega$ (and of course it
is infinite, otherwise trivial).
 
We define  a metric $\metricD$ on $G$ by 
\[\metricD(f,g)=\inf\{2^{-n}:f(n)\neq g(n)\mbox{ or } f^{-1}(n)\neq
g^{-1}(n)\}.\] 
Then:
\begin{enumerate}
\item $G$ is a complete separable metric space under $\metricD$, in fact a
topological group.  
\item If ${\mathbf d}$ is an $\omega$--sequence of members of $G\setminus\{ 
e_{G}\}$ converging to $e_{G}$, then for some (strictly increasing)
$\omega$--sequence ${\mathbf j}$ of natural numbers the pair $({\mathbf d},
{\mathbf j})$ satisfies   
\begin{enumerate}
\item[(*)] for any sequence $\langle w_{n}(x_1,x_2,\dots,x_{\ell_{1,n}};
y_1,y_2,\dots,y_{\ell_{2,n}}):n<\omega\rangle$ (i.e., an $\omega$--sequence
of group words) obeying ${\mathbf j}$ (see below) we can find a sequence
${\mathbf b}$ from $G$, that is $b_{n}\in G$ for $n<\omega$ such that  
\[b_{n}=w_{n}(d_{n+1},d_{n+2},\dots,d_{n+\ell_{1,n}};b_{n+1},b_{n+2},\dots,
b_{n+\ell_{2,n}})\qquad\mbox{ for any }n.\]  
\begin{enumerate}
\item[$(*)_1$] We say that $\langle w_{n}(x_1,x_2,\dots,x_{\ell_{1,n}};
y_1,y_2,\dots, y_{\ell_{2,n}}):n<\omega\rangle$ obeys ${\mathbf j}$ if: 
\item[$(*)_2$] for any $n^*,m^*<\omega$ we can find $i(0),i(1)$ such that  
\item[$(*)_3$] $m^*<i(0)$, $n^*<i(1)$, $i(0)<i(1)$, and $w_t$ is trivial
(which means $w_t=x_1$) for $t=j_{i(0)},j_{i(0)}+1,\dots,j_{i(1)}$, and
\[\sum\limits_{i=n,\dots,j_{i(0)}}\length(w_{i})<i(1)-i(0),\]  
where  
\item[$(*)_4$] the length of a word $w=w(z_1,\dots,z_{r})$ which in
canonical form is $z_{\pi(1)}^{t(1)} z_{\pi(2)}^{t(2)}\dots
z_{\pi(s)}^{t(s)}$ with $\pi$ a function from $\{1,2,\dots,s\}$ to
$\{1,\dots,r\}$  and $t(i)\in {\mathbb Z}$, is $\sum\limits_{i=1,\dots,s}
\vert t(i)\vert$
\end{enumerate}
\end{enumerate}
\end{enumerate}
\end{proposition}

\begin{proof}
(1)\quad Should be clear 
  
\noindent (2)\quad So we are given the sequence ${\mathbf d}$. We choose the
increasing sequence ${\mathbf j}$ of natural numbers by letting $j_0=0$,
$j_{n+1}$ be the first $j>j_{n}$ such that

\begin{enumerate}
\item[$(*)_5$] for every $\ell\le n$ and $m<j_{n}$ we have $d_\ell(m)<j$ and
$(d_\ell)^{-1}(m)<j$, and $[k\ge j\ \Rightarrow\ d_{k}(m)=m$], and $w_\ell$
mentions only $x_{r}, y_{r}$ with $r<j-j_{n}$.  
\end{enumerate}
Note that $j_{n+1}$ is well defined as the sequence ${\mathbf d}$ converges
to $e_{G}$, so for each $m<\omega$ for every large enough $k<\omega$ we have
$d_{k}(m)=m$.  
  
We shall prove that ${\mathbf j}=\langle j_n:n<\omega\rangle$ is as required
in part (2) of the proposition. So let a sequence $\langle w_{n}(x_1,x_2,
\dots,x_{\ell_{1,n}}; y_1,y_2,\dots,y_{\ell_{2,n}}):n<\omega\rangle$ of
group words obeying ${\mathbf j}$ be given (see $(*)_2 + (*)_3$ above). 
   
For each $k<\omega$ we define the sequence $\langle b^{k}_{n}:n<\omega
\rangle$ of members of $G$ as follows. For $n>k$ we let $b^{k}_n$ be $e_G$
and now we define $b^{k}_{n}$ by downward induction on $n\le k$ letting  

\begin{enumerate}
\item[$(*)_6$]  $b^{k}_{n}=w_{n}(d_{n+1},d_{n+2},\dots,d_{n+\ell_{1,n}};
b^{k}_{n+1},b^2_{n+2},\dots,b^2_{n+\ell_{2,n}})$.    
\end{enumerate}
Now  shall work on proving 

\begin{enumerate}
\item[$(*)_7$] for each $n^*,m^*<\omega$ the sequence $\langle b^{k}_{n^*}(
m^*):k<\omega\rangle$ is eventually constant.
\end{enumerate}
Why does $(*)_7$ hold? By the definition of obeying we can find $i(0),i(1)$
such that  

\begin{enumerate}  
\item[$(*)_8$] $m^*<i(0)$, $n^*<i(1)$, $i(0)<i(1)$, and $w_{t}$ is trivial
for $t=j_{i(0)},j_{i(0)}+1,\dots,j_{i(1)}$, and 
\[\sum\limits_{i=n,\dots,j_{i(0)}}\length(w_{i})<i(1)-i(0).\]
\end{enumerate}
For $n\in [n^*,j_{i(0)}]$ let 
\[t(n)=\sum_{i=n^*,\dots,n-1}\length(w_{i}).\] 
Now let $k(*)=^{\rm df} j_{i(1)+1}$; we claim that: 

\begin{enumerate}  
\item[$(*)_9$] if $k\ge k(*)$ and $s\ge j_{i(1)}$, then $b^{k}_{s}$
restricted to the interval $[0,j_{i(1)-1})$ is the identity.  
\end{enumerate}
[Why? If $s>k$, this holds by the choice of the $b^{k}_{s}$ as the identity
everywhere.  Now  we prove $(*)_9$ by downward induction on $s\le k$ (but of
course $s\ge j_{i(1)}$). But by the definition of composition of
permutations it suffices to show 

\begin{enumerate}  
\item[$(*)_{9a}$] every permutation mentioned in the word 
\[w_{s}(d_{s+1},\dots,d_{s+\ell_{1,s}},b^{k}_{s+1},\dots, b^{k}_{s+
\ell_{2,s}})\]
maps every $m<j_{i(1)-1}$ to itself. 
\end{enumerate}
Let us check this criterion. The $d_{s+\ell}$ for $\ell=1,\dots,\ell_{1,s}$
satisfies this as the indexes are $\ge j_{i(1)}$ and $m<j_{i(1)-1}$; now
apply the choice of $j_{i(1)}$.  
  
The $b^{k}_{s+1},\dots,b^{k}_{s+\ell_{2,s}}$ satisfy this by the induction
hypothesis on $s$. So the demands in $(*)_{9a}$ holds, hence we complete the
downward induction on $s$. So $(*)_9$ holds.] 

\begin{enumerate}  
\item[$(*)_{10}$] If $k\ge k(*)$ and $s\in [j_{i(0)},j_{i(1)}]$, then
$b^{k}_{s}$ is the identity on the interval $[0,j_{i(1)-1})$.  
\end{enumerate}
[Why? We prove this by downward induction; for $s=j_{i(1)}$ this holds by
$(*)_9$, if it holds for $s+1$, recall that $w_{s}$ is trivial, so
$b^{k}_{s}=b^{k}_{s+1}$, so this follows.] 

\begin{enumerate}
\item[$(*)_{11}$] For every $k\ge k(*)$ we have: for every $s\geq j_{i(0)}$,
the functions $b^{k}_{s},b^{k(*)}_{s}$ agree on the interval
$[0,j_{i(1)-1})$, and also $(b^k_s)^{-1}, (b^{k(*)}_s)^{-1}$ agree on this
interval. 
\end{enumerate}
[Why? For $s\geq j_{i(1)}$ by $(*)_9$; for $s\in [j_{i(0)},j_{i(1)})$ by
$(*)_{10}$.]  
  
\begin{enumerate}  
\item[$(*)_{12}$] For any $s\in [n^*,\omega)$ and $m<j_{i(1)-1}$ such that
\[[s<j_{i(0)}\ \Rightarrow\ m<j_{i(0)+t(s)}]\]
we have:  
\[\begin{array}{l}
k\ge k(*) (=j_{i(1)+1})\quad\mbox{ implies:}\\
b^{k}_{s}(m)=b^{k(*)}_{s}(m)\ \mbox{ and }\ (b^{k}_{s})^{-1}(m)=(
b^{k(*)}_{s})^{-1}(m).
  \end{array}\] 
\end{enumerate}
  
\noindent {\sc Case 1}:\quad $s$ is $\ge j_{i(0)}$.
 
\noindent [Why? This holds by $(*)_{11}$.] 

We prove this by downward induction on $s$ (for all $m$ and $k$ as there).  
  
\noindent {\sc Case 2}:\quad Proving for $s<j_{i(0)}$, assuming we have it
for all relevant $s'>s$ (and $s\ge n^*$ of course).  

\noindent Let $k\ge k(*)$ and we concentrate on proving $b^{k}_{s}(m)=
b^{k(*)}_{s}(m)$ as the proof of $(b^k_s)^{-1}(m)=(b^{k(*)}_s)^{-1}(m)$ is
the same. So   
\[b^{k}_{s}(m)=w_{s}(d_{s+1},\dots,d_{s+\ell_{1,s}},b^{k}_{s+1},\dots,
b^{k}_{s+\ell_{2,s}}).\]
So let us write this group expression as the product $u^k_{s,1}\dots u^k_{s, 
\length(w_s)}$, where each $u^k_{s,r}$ is one of $\{d_{s+1},\dots,d_{s+
\ell_{1,s}},b^{k}_{s+1},\dots,b^{k}_{s+\ell_{2,s}}\}$, or is an inverse of
one of them. 

For $r=0,1,\dots,\length(w_s)$ let $v^{k}_{s,r}=u^k_{s,1}\dots u^k_{s,r}$,
so $v^{k}_{s,r}\in G$ is the identity permutation for $r=0$ and is $w_{s}(
d_{s+1},\dots,d_{s+\ell_{1,s}},b^{k}_{s+1},\dots,b^{k}_{\ell_{2,s}})$ for
$r=\length(w_s)$. Hence it suffices to prove the following 
\begin{enumerate}  
\item[$(*)_{12a}$] if $r\in\{0,\dots,\length(w_{n})\}$ and $m<j_{i(0)+t(s)
+r}$, then $v^{k}_{s,r}(m)=v^{k(*)}_{s,r}(m)$.
\end{enumerate}
[Why does $(*)_{12a}$ hold? We do it by induction on $r$; now for $r=0$ the
permutation is the identity so trivial. For $r+1$ just note that because
each $u^k_{s,r}$ can map any $m<j_{i(0)+t(n)+r}$ only to numbers $m'<j_{i(0)
+t(n)+r+1}$ and that $(*)_{12}$ has been proved for $b^{k}_{s'},
b^{k(*)}_{s'}$ when $s'>s$ is appropriate.]  
  
So we have proved $(*)_{12a}$, and hence $(*)_{12}$ and thus also $(*)_7$. 

Lastly 
\begin{enumerate}  
\item[$(*)_{13}$] for each $n^*$ and $m^*<\omega$ the sequence $\langle
(b^{k}_{n^*})^{-1}(m^*):k<\omega\rangle$ is eventually constant.
\end{enumerate}
Why? Same as the proof of $(*)_7$.  

Together, we can defined for any $m,n<\omega$ the natural number
$b^*_{n}(m)$ as the eventual value of $\langle b^{k}_{n}(m):k<\omega
\rangle$. So $b^*_{n}$ is a well defined function from the natural numbers
to themselves (by $(*)_7$), in fact it is one-to-one (as each $b^{k}_{n}$
is) and is onto (by $(*)_{13}$), so it is a permutation of ${\mathbb
A}$. Clearly the sequence $\langle b^{k}_{n}:k<\omega\rangle$ converges to
$b^*_{n}$ as a permutation, the metric is actually defined on the group of
permutations of the family of members of ${\mathbb A}$, and $G$ is a closed 
subgroup; so $b^*_{n}$ actually is an automorphism of ${\mathbb A}$. 
Similarly the required equations   
\[b^*_{n}=w_{n}(d_{n+1},\dots,d_{n+\ell_{1,n}},b^*_{n+1},\dots,b^*_{n+
\ell_{2,n}})\]
hold. 
\end{proof}

\begin{proposition}  
\label{prop4}
The conclusion of Proposition 3 fails for any uncountable free group $G$. 
\end{proposition}
 
\begin{proof}
So let $Y$ be a basis of $G$ and as $G$ is a separable metric space there is
a sequence $\langle c_{n}:n<\omega\rangle$ of (pairwise distinct) members of
$G$ with $\metricD(c_{n},c_{n+1})<2^{-n}$. Let $d_{n}=(c_{2n})^{-1}
c_{2n+1}$, so $\langle d_{n} : n<\omega\rangle$ converges to $e_{G}$ and
$d_n\neq e_G$. Assume $\langle j_{n}:n<\omega\rangle$ is as in the
conclusion of Proposition \ref{prop3}, and we shall eventually get a
contradiction. Let $H$ be a subgroup of $G$ generated by some countable
$Z\subseteq Y$ and including $\{c_{n}:n<\omega\}$ and let $\langle a_{n}:n<
\omega\rangle$ list the members of $H$.  
  
Now 
\begin{enumerate}  
\item[$(*)_1$] $\langle d_{n}:n<\omega\rangle$ satisfies the condition also
in $H$. 
\end{enumerate}
[Why? As there is a projection from $G$ onto $H$ and $d_{n}\in H$.] 

For each $\nu\in {}^{\omega}\omega$ let ${\mathbf w}_{\nu}=\langle
w^{\nu}_{n}:n<\omega\rangle$, where $w^{\nu}_{n}=w^{\nu}_{n}(x_1,y_1)=x_1
(y_1 )^{\nu(n)}$, so this is a sequence of words as mentioned in Proposition
\ref{prop3}.  

\begin{enumerate}
\item[$(*)_2$] The set of $\nu\in {}^{\omega}\omega$ which obey ${\mathbf
j}$ is co-meagre. 
\end{enumerate}
[Why? Easy; for each $n^*,m^*<\omega$ the set of ${\mathbf w}$'s which fail
the demand for $n^*,m^*$ is nowhere dense (and closed), hence the set of
those failing it is the union of countably many nowhere dense sets, hence is
meagre.]  
  
\begin{enumerate}
\item[$(*)_3$] For each $a\in H$ the family of $\nu\in {}^{\omega}\omega$
such that there is solution ${\mathbf b}$ for $({\mathbf d},{\mathbf w})$
in $H$ satisfying $b_0=a$ is nowhere dense. 
\end{enumerate}
[Why? Given a finite sequence $\nu$ of natural numbers note that for any
sequence $\rho\in {}^{\omega}\omega$ of which $\nu$ is an initial segment
and solution ${\mathbf b}$ satisfying $b_0=a$, we can show by induction
on $n\le\lh(\nu)$ that $b_{n}$ is uniquely determined, call it
$b[n,\nu,{\mathbf d}]$. Now, if $b[\lh(\nu),\nu,{\mathbf d}]$, which is a
member of $G$, is not $e_{G}$, then for some $m<\omega$ it has no $t$-th
root and we let $\nu_1=\nu\conc \langle t\rangle$ and we are done. If not,
letting $\nu_0=\nu\conc\langle 0 \rangle$, also $b[\lh(\nu)+1,\nu_1,{\mathbf
d}]$ is well defined and equal to $d_{\lh(\nu)}$, hence is not $e_{G}$. 
Hence for some $t<\omega$ has no $t$-th root, so $\nu_1=\nu_0\conc\langle
t\rangle$ is as required.]  
  
Now we can finish the proof of Proposition \ref{prop4}: just by $(*)_2+(*)_3
+$ Baire Theorem, for some $\nu\in {}^\omega\omega$, the sequence ${\mathbf
w}$ of group words obeying ${\mathbf j}$, there is no solution in $H$, hence
no solution in $G.$  
\end{proof}

\begin{proof}[Proof of Theorem \ref{Thm1}] 
Follows by Propositions \ref{prop3}, \ref{prop4}.
\end{proof}

\begin{conrem}
\begin{enumerate}
\item[(A)] In the proof of Proposition \ref{prop4} we do not use all the
strength of ``$G$ is free''. E.g., it is enough to assume:  
\begin{enumerate}
\item[(a)] if  $g\in G$, $g\neq e_{G}$, then for some $t>1$, $g$ has no
$t$-th root (in $G$),  
\item[(b)] if $X$ is a countable subset of $G$, then there is a countable
subgroup $H$ of $G$ which includes $X$ and there is a projection from $G$
onto $H$,  
\item[(c)] $G$ is uncountable. 
\end{enumerate}
The uncountable free abelian group fall under this criterion; in fact by
Proposition \ref{prop3}, $G$ is ``large'', ``rich''.
\item[(B)] What about uncountable structures? Sometimes a parallel result
holds: if $\lambda=\beth_{\omega}$, replacing countable by ``of cardinality
$\le\beth_{\omega}$''. More generally, assume $\aleph_0<\lambda=
\sum\limits_{n<\omega}\lambda_n$ and 2$^{\lambda_n}<2^{\lambda_{n+1}}$ for
$n<\omega$, hence $\mu=^{\rm df}\sum\limits_{n<\omega} 2^{\lambda_n}<
2^\lambda$; and we have 
\begin{enumerate}
\item[$(*)$] if ${\mathbb A}$ is a structure with exactly $\lambda$ elements
and $G$ is its group of automorphisms, then $G$ cannot be a free group of
cardinality $>\mu$.  
\end{enumerate}
The proof is similar, but now w.l.o.g the set of elements of ${\mathbb A}$
is $\lambda=\{\alpha:\alpha<\lambda\}$ and we define $\metricD$ by  
\[\begin{array}{ll}
\metricD(f,g)=\inf\{2^{-n}:&\mbox{there is }\alpha<\lambda_n\mbox{ such
that}\\
&\mbox{for some }(f',g')\in\{(f,g),(f^{-1},g^{-1}),(g,f),(g^{-1},f^{-1})\}\\
&\mbox{one of the following possibilities holds}\\
&\mbox{(a)\ \ for some }m<\omega\mbox{ we have }f'(m)<\lambda_m\le g'(m),\\   
&\mbox{(b)\ \ }f'(n)< g'(n)<\lambda_n\}.
  \end{array}\]
Under this metric, $G$ is a complete metric space with density
$\le\sum\limits_{n<\omega} 2^{\lambda_n}=\mu$, and the conclusion of
Proposition \ref{prop3} holds. 
\item[(C)] By \cite{JShT:654}, for $\kappa=\kappa^{<\kappa}>\aleph_0$ there
is a forcing adding such a group and not changing cardinalities or
cofinalities. For ZFC result see \cite{Sh:F442}. 
\end{enumerate}
\end{conrem}


\end{document}